\documentclass[12pt]{amsart}

\usepackage{color}
\usepackage[centertags]{amsmath}
\usepackage{amsfonts}
\usepackage{amsthm}
\usepackage{newlfont}
\usepackage{amscd}
\usepackage{amsgen}
\usepackage{amssymb}
\usepackage{stmaryrd}
\usepackage{amssymb,amsmath}
\usepackage{amscd}
\usepackage{enumerate}
\usepackage{color}
\usepackage{xypic}
\usepackage{graphicx}

\newlength{\defbaselineskip} 
\theoremstyle{plain}
\newtheorem{thm}{Theorem}[section]
\newtheorem{cor}[thm]{Corollary}
\newtheorem{con}[thm]{Conjecture}
\newtheorem{df}[thm]{Definition}
\newtheorem{lema}[thm]{Lemma}
\newtheorem{obs}[thm]{Proposition}
\newtheorem{exm}[thm]{Example}

\newtheorem{fact}[thm]{Fact}

\newtheorem{rem}[thm]{Remark}
\newtheorem{pr}{Program}
\theoremstyle{definition} 
\theoremstyle{definition} \newtheorem{ex}{Example}[section] %

 \numberwithin{equation}{section}
\def\p{\mathbb{P}}

\def\z{\mathbb{Z}}
\def\n{\mathbb{N}}
\def\Z{\mathbb{Z}}

\def\c{\mathbb{C}}

\DeclareMathOperator{\Imi}{Im}
\DeclareMathOperator{\End}{End}

\def\p{\mathbb{P}}
\def\a{\mathbb{A}}
\def\ob{\begin{obs}}
\def\kob{\end{obs}}
\def\dow{\begin{proof}}
\def\kdow{\end{proof}}

\def\tw{\begin{thm}}
\def\ktw{\end{thm}}
\def\hip{\begin{con}}
\def\khip{\end{con}}
\def\fa{\begin{fact}}
\def\kfa{\end{fact}}
\def\lem{\begin{lema}}
\def\klem{\end{lema}}
\def\ex{\begin{exm}}
\def\prog{\begin{pr}}
\def\kprog{\end{pr}}
\def\wn{\begin{cor}}
\def\kwn{\end{cor}}
\def\uwa{\begin{rem}}
\def\kuwa{\end{rem}}
\def\kex{\end{exm}}
\def\dfi{\begin{df}}
\def\kdfi{\end{df}}
\definecolor{zielony}{rgb}{0.5, 0.9, 0.1}
\definecolor{czerwony}{rgb}{0.9, 0.2, 0.1}
\definecolor{niebieski}{rgb}{0.3, 0.1, 0.9}

\begin{document}
\title{Toric geometry of the 3-Kimura model for any tree}
\author{Mateusz Micha\l ek}
\thanks{The author is supported by a grant of Polish MNiSzW (N N201 413 539).}
\begin{abstract}
In this paper we present geometric features of group based phylogenetic models. We address a long standing problem of determining the ideal of the claw tree \cite{SS}, \cite{DK}. We focus on the 3-Kimura model. In particular we present a precise geometric description of the variety associated to any tree on a Zariski open set. This set contains all biologically meaningful points.
The result confirms the conjecture of Sturmfels and Sullivant \cite{SS} on the degree in which the ideal associated to 3-Kimura model is generated on that set.
\end{abstract}
\maketitle
\section{Introduction}
To a tree $T$ and an abelian group $G$ one can associate an algebraic variety $X(T,G)$. The construction of this variety is motivated by biology. The tree $T$ describes the evolution of species and the group $G$ distinguishes a model of evolution. By a model of evolution we mean certain constraints on the probabilities of mutation. In the article we do not consider the details of biological aspects of this construction. They are well described in \cite{PS}, or in a short paper \cite{4aut}. Our main object of study is the variety $X(T,G)$ and its geometry. Nevertheless some of the applications of the results will be suggested. Moreover the case of the group $G=\z_2\times\z_2$ will be of particular interest. This group corresponds to the 3-Kimura model appearing in theoretical biology.

From the point of view of algebraic geometry the variety $X(T,G)$ has got a lot of interesting properties. The first results are presented in \cite{ES} and \cite{SSE}. In the modern mathematical language, the authors prove that $X(T,G)$ is toric. We want to stress that throughout the article we do not assume that a toric variety is normal. Directly speaking $X(T,G)$ is the closure of the image of an algebraic map given by monomials.

One of the main advantages of toric geometry is its interaction with discrete geometry. The use of lattices, polytopes and fans allows to precisely describe toric varieties and do exact computations \cite{Ful}, \cite{Cox}. The varieties $X(T,G)$ also have a lot of connections with discrete objects - trees, groups, polytopes. Hence, the study of them often involves combinatorics, toric geometry or even representation theory. One of the examples is presented in \cite{BW} where the binary model $G=\z_2$ is studied. The major result presented there says that if two trees $T_1$ and $T_2$ have the same number of leaves, then the varieties $X(T_1,\z_2)$ and $X(T_2,\Z_2)$ belong to the same flat family. Unfortunately this is a special feature of the group $\z_2$ \cite{Kaie}, \cite[Section 3.2]{DBM}. One can also find deep connections with the representation theory. The algebras of the considered varieties are connected with algebras of conformal blocks \cite{Man}. Further interesting results can be found in \cite{Xu}, \cite{graphs}, \cite{AllRhMarkov}, \cite{DK} and \cite{sull}.

The 3-Kimura model is much more complicated then the binary model. In particular the Hilbert polynomial of the associated variety depends on the topology of the tree, not only the number of leaves. Still this model is one of the simplest that are of interest to biologists.

The main aim of this article is to present geometric constructions that lead to a better understanding of the variety $X(T,\z_2\times\z_2)$. One of the motivations is the conjecture of Sturmfels and Sullivant. It claims that the ideal of $X(T,\z_2\times\z_2)$ is always generated in degree $4$ \cite[Conjecture 2]{SS}. This conjecture seems more algebraic than geometric. However it can be presented in a purely geometric language in terms of scheme-theoretic intersection of varieties. Let $K_{1,l}$ be the claw tree with one inner vertex and $l$ leaves.
\hip[Conjecture 4.6 \cite{DBM}, cf.~Conjecture 2 \cite{SS}]\label{hipws}
The variety $X(K_{1,l},\z_2\times\z_2)$ is a scheme-theoretic intersection of the varieties $X(T,\z_2\times\z_2)$, where $T$ is any tree with $l$ leaves different from $K_{1,l}$.
\khip
The above conjecture can be stated for any phylogenetic model. It could provide the description of the ideal of claw trees for many relevant models. The reduction to claw trees for the general Markov model was presented in \cite{AllRhMarkov} and for equivariant models in \cite{DK}. The application of the results for 3-Kimura model is presented in Section \ref{Ap}. To stress the importance of the conjecture we cite \cite{DK}:

"We have now reduced the ideals of our equivariant models to those for stars, and
argued their relevance for statistical applications. The main missing ingredients
for successful applications are equations for star models. These are very hard to
come by (...)".

Let us fix a coordinate system on $\a^n\supset X(T,\z_2\times\z_2)$.
Let $B\subset\a^n$ be the Zariski open subset containing points with all coordinates different from zero. In this paper we prove that Conjecture \ref{hipws} holds on $B$. It follows that Conjecture 2 of \cite{SS} holds after a localization of the algebra with respect to all the coordinates. In particular all of the previous conjectures hold if and only if the intersection of the varieties $X(T,\z_2\times\z_2)$ is reduced and irreducible. The main idea of the proof is to transfer the properties of $X(T,\z_2)$ to $X(T,\z_2\times\z_2)$. This is done by an interplay between toric geometry and combinatorics. The best reference for such techniques is \cite{Stks}.

Among the applications we give all biologically meaningful invariants for the 3-Kimura model on any tree. We also discuss the identifiability issue. The main aim of the paper is to interest toric geometers and algebraic combinatorialists in the arising problems. Hence, we postpone the discussion of applications to Section \ref{Ap}.

Let us briefly discuss the structure of the paper. Section \ref{Not} settles the notation for the rest of the paper. All introduced notions are formal mathematical objects and much of the biological motivation is skipped. The most important are Definitions \ref{Socket} and \ref{Network} that are far from being standard in algebraic phylogenetics. Therefore an expert in the discipline is advised at least to take a look at these two definitions. People with mathematical background and no experience in biology should be able to follow this section\footnote{Although we admit that some of the constructions may appear unnatural.} and are encouraged to do so.

In Section \ref{BProp} we introduce objects typical of toric geometry. Some of the presented facts may be well known to toric geometers. We include them as an introduction to the proof of the main theorem. We also explain special features of varieties $X(T,G)$ and interactions between them.

Section \ref{Cat} describes relations between groups $G$ and varieties $X(T,G)$. We try to avoid general theory that will not be needed. Instead, we introduce basic concepts that are later used directly to prove the main theorem. We are interested in the interactions between morphisms of groups and morphisms of varieties. We are not setting a general categorical environment for phylogenetic models. This can be done, but involves more sophisticated language and is not required for the rest of the paper. The restriction to the case where $G$ is abelian allows to obtain short, natural statements.

In Section \ref{MT} we bring all the ideas together to prove the main theorem. The following Section \ref{Ap} discusses some of the possible applications. In fact, the dense open subset $B$ contains all points interesting form the biological point of view \cite{CFS}. This section should be particulary important to people interested in algebraic phylogenetics. It is self contained and results presented there are not used in the rest of the paper. However in this section we allow ourselves to use standard notions from phylogenetics without introducing them.

The last Section \ref{Sc} describes some partial results concerning the complement of $B$. It is quite technical. It contains a few results in favor of Conjecture \ref{hipws}. In Appendix we present two results about the binary model.
\section*{Acknowledgements}
I would like to thank very much Jaros\l aw Wi\'sniewski for introducing me to the subject and guiding me through it. I would also like to thank Marta Casanellas for important remarks.

\section{Notation}\label{Not}
We will be dealing with algebraic varieties associated to phylogenetic models. These varieties are always given as closures of the image of a parametrization map. A short introduction to the topic can be found in \cite{4aut}. In this paper we will be interested in the geometry of algebraic varieties. Hence we omit their detailed construction motivated by biology. It can be found for example in \cite[Section 2]{mateusz}, \cite{SS} or in a nice theoretical setting in \cite{DK}. We briefly present a purely algebraic construction of varieties associated to $G$-models. The proofs of these results can be found in \cite{mateusz}. We present them to make the paper self-contained for any algebraic geometer with a limited or no knowledge of biology. We hope that arising geometric problems will interest toric geometers, even if they are not dealing with phylogenetics.

We will be defining objects that will depend on a tree $T$ and a group $G$. For any object $O$ if we want to stress its dependence on either $T$ or $G$ we write them in the upper index: $O^{T,G}$.
\subsection{A variety associated to a $G$-model}
Let $S$ be a finite set.
\dfi[Space $W$]\label{W}
 We define $W$ to be a complex vector space spanned freely by elements of $S$. More precisely
$W=\oplus_{a\in S} \c_a,$
where $\c_a$ is a field of complex numbers corresponding to one dimensional vector space spanned by $a\in S$.
\kdfi
Suppose that a group $G$ acts on the set $S$, hence also on the vector space $W$. This induces a $G$ action on $\End (W)$. For $f\in\End(W)$, $g\in G$ and $x\in W$ we have $(gf)(x)=g(f(g^{-1}x))$.
\dfi[Space $\widehat W$]
We define $\widehat W\subset \End(W)$ as the maximal space of endomorphisms invariant with respect to the $G$ action.
\kdfi
Let us consider the following setting. Assume that $G$ contains a normal, abelian subgroup $H$ that acts transitively and freely on the set $S$. This is a technical assumption. It is crucial however to prove that the varieties we obtain are toric. Models obtained in this way are called $G$-models \cite{mateusz}, \cite{WDB}. The majority of the paper concerns the binary model and the 3-Kimura model. They both correspond to the case when $G$ is abelian.
\dfi[General group based model]
We will say that we are dealing with a general group-based model if $G$ is abelian. Moreover we assume that $S$ is a torsor for $G$. In this setting the vector space $W$ is the regular representation of $G$.
\kdfi
\uwa\emph{
Models coming from biology do not have to be given by a group action. However the class of $G$-models is the largest class known to the author that gives rise to toric varieties.
}\kuwa
In the following discussion we assume that we are dealing with a general group-based model. Using more complicated language similar results can be obtained for $G$-models \cite{mateusz}. The restriction to abelian groups in particular allows us to work with group elements, instead of sets of characters.
Suppose that we are given a rooted tree $T$ with edges directed from the root.
\dfi[Sets $L$, $V$, $N$ and $E$]
Let $L$, $V$, $N$ and $E$ be respectively the set of leaves, vertices, nodes and edges of the tree $T$. We have $V=L\cup N$. We identify leaves with edges adjacent to them.
\kdfi
\dfi[Socket]\label{Socket}
A socket is a function $s:L\rightarrow G$ such that $\sum_{l\in L} s(l)=0$, where $0$ is the neutral element of the group $G$.
\kdfi
\ex\label{socket}\emph{
Let us consider the group $G=\z_3$ and the following tree:
$$
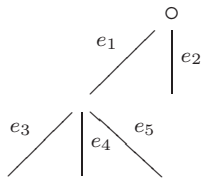
\xymatrix{
&&\circ\ar@{-}[dl]_{e_1}\ar@{-}[d]^{e_2}\\
&\ar@{-}[dl]_{e_3}\ar@{-}[d]^{e_4}\ar@{-}[dr]^{e_5}&\\
&&\\
}$$
Here $e_2$, $e_3$, $e_4$ and $e_5$ are leaves. An example of a socket is an association $e_2\rightarrow 1$, $e_3\rightarrow 1$, $e_4\rightarrow 2$, $e_5\rightarrow 2$.}
\kex
Suppose that we have associated group elements to edges adjacent to a node $n\in N$. We say that a signed sum around this node is trivial if the element associated to the incoming edge equals the sum of elements associated to outcoming edges. For the root we want the sum of elements associated to outcoming edges to be the neutral element.
\dfi[Network]\label{Network}
A network is a function $f:E\rightarrow G$ such that the signed sum of associated elements around each node is trivial.
\kdfi
\ex\emph{
We consider the same tree as in Example \ref{socket}. We can make a network using the same association and extending it by $e_1\rightarrow 2$.}
\kex
\ob\label{bijekcja}
There is a natural bijection between sockets and networks.
\kob
\dow
Let $n$ be a network. We know that the signed sum $S(v)$ around each inner vertex $v$ is the neutral element element. Hence $\sum_{v\in N} S(v)=e$, where $e$ is the neutral element. Let us consider an edge directed from $v_1$ to $v_2$, where $v_1,v_2\in N$. Let us note that the group elements $n(v_1,v_2)$ and $n(v_1,v_2)^{-1}$ appear in $S(v_1)$ and $S(v_2)$. We see that $\sum_{v\in N} S(v)=\sum_{l\in L} n(l)$. This means that a restriction of the network to leaves gives a socket.

Given a socket $s$ we can define a function $n:E\rightarrow G$ using the condition of summing up to the neutral element around inner edges. The only nontrivial thing is to notice that the sum around the root also gives the neutral element. This follows from the previous equality $\sum_{v\in N} S(v)=\sum_{l\in L} n(l)$ and the fact that $S(v)=e$ for each node $v$ different from the root.
\kdow
To each edge $e\in E$ we associate a vector space $\widehat W_e\cong \widehat W$. It is an easy exercise that $\widehat W_e$ is $|G|$-dimensional. We assume that the coordinates of $\widehat W_e$ are indexed by group elements. In the same way to each leaf $l$ we associate a vector space $W_l\cong W$ with coordinates indexed by elements of $G$.
\dfi[Spaces $\widehat W_E$ and $W_L$] We define the following spaces:
$$W_L=\bigotimes_{l\in L} W_l,\qquad \widehat W_E=\bigotimes_{e\in E}\widehat W_e.$$
We call $W_L$ the space of states of leaves and $\widehat W_E$ the parameter space.
\kdfi
An element of a base of $W_L$ can be given by a tensor product of base elements of each $W_l$. The base of $W_l$ corresponds to group elements. Hence elements of the base of $W_L$ correspond to associations of group elements to leaves. In the same way elements of the base of $\widehat W_E$ are represented by associations of group elements to edges. We define two subspaces spanned by subsets of distinguished basis.
\dfi[Spaces $\widetilde W_E$, $\widetilde W_L$]
The subspace $\widetilde W_E\subset \widehat W_E$ is spanned by basis elements corresponding to associations that form a network.
The subspace $\widetilde W_L\subset W_L$ is spanned by basis elements corresponding to associations that form a socket.
\kdfi
Thanks to Proposition \ref{bijekcja} we have got a natural isomorphism $\widetilde W_E\cong\widetilde W_L$.
We have got the following map:
$$\psi: \prod_{e\in E} W_e\rightarrow \widehat W_E\rightarrow \widetilde W_E\cong\widetilde W_L.$$
The first arrow is the tensor product and the second is the projection onto the subspace.
The induced projective maps are the following:
$$\psi_{\p}:\prod_{e\in E}\p(W_e)\rightarrow \p(\widehat W_E)\dashrightarrow \p(W_L).$$
The first map is the Segre embedding and the second is a projection. Note that it follows from definition that the maps $\psi$ and $\psi_{\p}$ are given by monomials in the given basis.
The closure of the image of $\psi_{\p}$ (resp.~$\psi$) is the projective (resp.~affine) variety associated to the model. This is the main object of our study. We will denote it by $\p(X(T,G))$ (resp.~$X(T,G)$).
\section{Basic properties of algebraic varieties associated to phylogenetic models}\label{BProp}
\subsection{Toric varieties}
First we have to point out that we do not assume that a toric variety has to be normal. We only assume that a torus acts on a variety and one of the orbits is dense. This setting is most common when dealing with applications. Much information can be found in \cite{Stks}. The main drawback of this approach is that the varieties we consider will not be given by a fan. However, still they can be represented by polytopes, that do not have to be normal. For this reason we will often work with the character lattice $M$ instead of the one parameter subgroup lattice $N$.

Each vector space can be considered as a toric variety. We fix a coordinate system. The set of points with nonzero coordinates is an algebraic torus. It acts on the vector space. Moreover the dense torus orbit can be identified with the torus. As each vector space we have considered so far had a distinguished basis we can associate to it the character lattice of a torus.
\dfi[Lattices $M_S$, $M_e$, $M_E$]
For en edge $e\in E$ the lattice $M_e$ is the character lattice of the torus acting on $\widehat W_e$.
The lattice $M_E$ is the character lattice of the torus acting on $\prod_{e\in E} W_e$.
The lattice $M_S$ is the character lattice of the torus acting on $\widetilde W_E\cong \widetilde W_L$.
\kdfi
Let us note that the coordinate system on the vector space distinguishes the basis of the lattice. The basis of each lattice $M_e$ is indexed by group elements. As $M_E=\bigotimes_{e\in E} M_e$ the basis of $M_E$ is indexed by pairs $(e,g)$ where $e$ is an edge and $g$ a group element. The basis elements of $M_S$ corresponds to sockets or networks. The morphism $\psi: \prod_{e\in E} W_e\rightarrow \widetilde W_E$ is a toric morphism.
\dfi[Morphism $\widehat\psi$]
The morphism $\widehat\psi: M_S\rightarrow M_E$ is the morphism of lattices induced by $\psi$. 
\kdfi
In this setting the description of $\widehat\psi$ is particulary simple. Let $f_n\in M_S$ be a basis vector corresponding to a network $n$. The element $\widehat\psi(f_n)$ will be an element of the unit cube in $M_E$. Let $h_{(e,g)}\in M_E$ be the basis vector indexed by a pair $(e,g)\in E\times G$ and let $h^*_{(e,g)}$ be its dual. We have:
$$
h^*_{(e,g)}(f_n)=
\begin{cases}
1\text{  if  }n(e)=g\\
0\text{  otherwise.}
\end{cases}
$$
We come to the \textbf{most important definition of this section}.
\dfi[Polytope $P$]
We define the polytope $P\subset M_E$ to be the image of the basis of $M_S$ by $\widehat\psi$. In other words the vertices of the polytope $P$ correspond to networks. More precisely each vertex has got $1$ on coordinates indexed by pairs that form a network and $0$ on other coordinates.
\kdfi
\ex\emph{
Let us consider the tree $T$ with one inner vertex and three leaves $l_1,$ $l_2$ and $l_3$. Let $G=\z_2$. The lattice $M_S$ is the 4 dimensional lattice generated freely by vectors $e_{(0,0,0)},$ $e_{(1,1,0)},$ $e_{(1,0,1)},$ $e_{(0,1,1)}$ that correspond to sockets/networks on $T$. The lattice $M_E$ is a 6 dimensional lattice with basis vectors $f_{(l_i,g)}$ with $1\leq i\leq 3$ and $g\in \z_2$. We have $\widehat\psi(e_{(a,b,c)})=f_{(l_1,a)}+f_{(l_2,b)}+f_{(l_3,c)}$. Hence each vertex of $P$ will have three coordinates equal to zero and three to one. Let us consider the base of $M_E$ in the following order $f_{(l_1,0)},f_{(l_1,1)},\dots,f_{(l_3,0)},f_{(l_3,1)}$. The vertex corresponding to $e_{(0,0,0)}$ is $(1,0,1,0,1,0)$. In the same order $e_{(1,1,0)}\rightarrow (0,1,0,1,1,0)$, $e_{(1,0,1)}\rightarrow(0,1,1,0,0,1)$ and $e_{(0,1,1)}\rightarrow (1,0,0,1,0,1)$. These are of course all vertices of $P$.}
\kex
The polytope $P$ is the polytope associated to the toric variety $X(T,G)$. The algebra of this variety is the algebra associated the semigroup generated by $P$ in $M_E$. Note that $P$ does not have to generate the lattice $M_E$.
\dfi[Lattice $\widehat M_E$]
We define the lattice $\widehat M_E$ as a sublattice of $M_E$ generated by $P$.
\kdfi
The lattices defined so far corresponded to affine objects. A rational map from a vector space to its projectivization is well defined on points with non zero coordinates. Hence it induces a surjective morphism of tori, what corresponds to an injective morphism of character lattices.
\dfi[Degree functions $\deg_e$]
Note that for a character lattice $M$ with a distinguished basis we can define a function $\deg:M\rightarrow \z$ that sums up coordinates. The degree of a lattice element is the degree of the monomial function associated to it.
For lattices $M_e$ the corresponding degree functions are denoted by $\deg_e$.
\kdfi
\dfi[Lattices $M_{S,0}$, $M_{E,0}$ and $\widehat M_{E,0}$]
For a lattice $M_S$ we define $M_{S,0}$ as a sublattice of elements with the sum of coordinates equal to zero.
In particular $M_{S,0}$ is the character lattice of the torus acting on $\p(\widetilde W_E)$.
We define $M_{E,0}$ as a sublattice of $M_E$ defined by equalities $\deg_e=0$ for each $e$. This is the character lattice of the torus acting on $\prod\p(W_e)$.
We also define $\widehat M_{E,0}:=M_{E,0}\cap \widehat M_E$. This is the character lattice of the torus acting on the projective toric variety $\p(X(T))$.
\kdfi
\subsection{Intersection of tori}\label{Intoftor}
To prove the main theorem we will also need to briefly review the theory of intersection of tori. Suppose that we have got tori $T_i$ that are subtori of a torus $T$. Let $M_i$ be the character lattice of the torus $T_i$ and $M$ the character lattice of the torus $T$. The inclusion $T_i \hookrightarrow T$ gives a surjective morphism $M\rightarrow M_i$. The kernel $K_i$ of this morphism corresponds to characters of $T$ that are trivial on $T_i$. In the algebra of $T$ the ideal of $T_i$ is generated by elements of the type $k-1$ for $k\in K_i$. Let $T'$ be the scheme theoretic intersection of the tori $T_i$. Let $M'$ be a sublattice of $M$ generated by the lattices $K_i$.
\fa\label{int}\emph{
The algebraic set $T'$ is irreducible iff $M'$ is saturated. The components of $T'$ are tori and the number of components equals the index of $M'$ in its saturation in $M$. If $M'$ is saturated then $M/M'$ is the character lattice of $T'$. The functions $k-1$ for $k\in M'$ generate the ideal of $T'$.
}\kfa
\subsection{Interactions between trees and varieties}
We can define an order on trees with $l$ leaves as follows. We say that $T_1\leq T_2$ if $T_1$ can be obtained form $T_2$ by a series of contractions of inner edges. Here by an edge contraction we mean identifying two vertices of a given edge. The smallest tree with $l$ leaves is the claw tree $K_{1,l}$ with one inner vertex. This is a part of a construction of the tree space \cite{sptr}.
We fix an abelian group $G$.
\ob\label{inkluzja}
If $T_1\leq T_2$ then $X(T_1,G)\subset X(T_2,G)$.
\kob
\dow
Although the statement is very easy we believe that the following discussion may be helpful to better understand the rest of the paper.
Both trees have got the same number of leaves, so we can make a natural bijection between their sockets. This gives an isomorphism of the ambient spaces $\widetilde W_E$.
As $T_1\leq T_2$ we can make an injection from the edges of $T_1$ to the edges of $T_2$. Note that a network on $T_2$, restricted to the edges of
$T_1$ is a network on $T_1$. This gives us a projection $\pi:M_{E}^{T_1} \twoheadrightarrow M_E^{T_2}$. The map $\pi$ simply forgets the coordinates indexed by $(e,g)$, where $e$ is an edge of $T_2$ not corresponding to an edge of $T_1$. Moreover the projection of $P^{T_2}$ is equal to $P^{T_1}$. The following diagram commutes:
$$
\xymatrixrowsep{5pt}
\xymatrix{
&M_E^{T_2}\ar@{->>}[dd]\\
M_S\ar[ru]\ar[dr]&\\
&M_E^{T_1}.\\
}
$$
Any relation between the vertices of $P^{T_2}$ is also a relation between the vertices of $P^{T_1}$. Hence any polynomial in the ideal of $X(T_2,G)$ is also in the ideal of $X(T_1,G)$.
\kdow
The surjective morphism of algebras corresponding to the inclusion of varieties is given by the restriction of the surjective morphism between $M_E^{T_2}$ and $M_E^{T_1}$ to the cones spanned by polytopes $P^{T_2}$ and $P^{T_1}$.

It is natural to ask what is the relation between $X(T_0,G)$ and the scheme theoretic intersection of all $X(T,G)$ for $T_0<T$. The first conjecture is that if there exists at least one $T>T_0$, then they are equal. An equivalent conjecture was stated in \cite{DBM}.
\hip\label{glhip}
The ideal of the variety associated to the claw tree $K_{1,l}$ is the sum of all the other ideals associated to trees with $l$ leaves.
\khip
So far we only know that the answer is positive for $G=\z_2$ \cite{Sonja}, \cite{SS}, \cite{DBM}. In the case particulary interesting for biologists $G=\z_2\times\z_2$ we know that this conjecture implies Conjecture 2 form \cite{SS} and is equivalent to it for $l> 8$. For other groups Conjecture \ref{glhip} is equivalent to Conjecture $1$ from \cite{SS} for $l$ large enough.

Conjecture \ref{glhip} can be stated for any phylogenetic model, not necessarily given by a group. In particular for a general Markov model.
One would be also interested to know exactly what is an intersection of a few varieties associated to different trees. In particular how many ideals do we have to sum to obtain the ideal associated to the claw tree. For the binary model some answers can be found in \cite{DBM}. A stronger result can be found in Appendix \ref{A}. One could also hope that the intersection of $X(T_1,G)$ and $X(T_2,G)$ is equal to $X(T,G)$ where $T$ is the largest tree smaller then $T_1$ and $T_2$. A counterexample to this statement is presented in Appendix \ref{B}.

\section{Category of general group based models}\label{Cat}
In this section we will establish connections between morphisms of groups and morphisms of corresponding varieties. Once again our main aim is application in geometry. We are building the set up of the next section. That is why we restrict to special cases. This reduces the complexity of the language but still gives a geometric insight. Let us fix a tree $T$.
\dfi[Groups $\mathfrak{S}$ and $\mathfrak{N}$]
By edgewise action sockets and networks on a given tree form a group. These groups will be denoted respectively by $\mathfrak{S}$ and $\mathfrak{N}$.
\kdfi
Let $f:G_1\rightarrow G_2$ be a morphism of abelian groups. It induces morphisms $\mathfrak{S}^{G_1}\rightarrow\mathfrak{S}^{G_2}$ and $\mathfrak{N}^{G_1}\rightarrow\mathfrak{N}^{G_2}$. This gives the following commutative diagram :
\begin{equation}\label{diag}\end{equation}
$$
\xymatrix{
M_S^{G_1}\ar[r]\ar[d]&M_E^{G_1}\ar[d]\\
M_S^{G_2}\ar[r]&M_E^{G_2}.\\
}
$$
Hence the morphism $M_E^{G_1}\rightarrow M_E^{G_2}$ of character lattices restricts to cones over polytopes. This gives a morphism of algebras of associated varieties. The morphism $M_S^{G_1}\rightarrow M_S^{G_2}$ restricts to positive quadrants of both lattices. Hence we get a morphism of ambient spaces $\widehat f:\widetilde W_L^{G_2}\rightarrow\widetilde W_L^{G_1}$ compatible with morphism of varieties $\widehat f:X(T,G_2)\rightarrow X(T,G_1)$. This gives a contravariant functor from the category of abelian groups to the category of embedded affine toric varieties. Moreover if $f$ is injective (resp.~sujective) then $\widehat f$ is dominant (resp.~injective). The second assertion is an easy exercise.

We also need the following setting. Suppose that we have morphisms $f_i:G\rightarrow G_i$ for $i=1,\dots,m$. Just as above this gives us a morphism of embedded varieties $\prod X(T,G_i)\rightarrow X(T,G)$. If the product $f_1\times\dots\times f_m$ is sujective then the morphism of varieties is injective. However in general if the product $f_1\times\dots\times f_m$ is injective the morphism of varieties does not have to be dominant.

\section{Main Theorem}\label{MT}
Our aim is to prove the following conjecture for $G=\z_2\times\z_2$.
\hip\label{drhip}
The dense torus orbit of the toric variety $X(K_{1,l},G)$ is an intersection of the dense torus orbits of varieties $X(T,G)$, where $T$ is any tree with $l$ leaves different form the claw tree.
\khip
Note that all dense torus orbits are contained in the dense torus orbit $O$ of the projective (or affine) ambient space. In the algebraic set $O$ all the considered orbits are closed subschemes. Hence Conjecture \ref{drhip} can be regarded in a set-theoretic or in a scheme-theoretic version. Both of them are equivalent. This follows for example from a more general statement \cite[Corollary 2.2]{EiS} and is particulary simple in toric case. However because the proofs of both versions are basically the same for $G=\z_2\times\z_2$ we have decided to include both. Moreover this also gives an idea how the elements of the ideal of $X(K_{1,l},\z_2\times\z_2)$ can be generated by elements of ideals of $X(T,\z_2\times\z_2)$.

The main idea of the proof is to extend the results known for binary models to the 3-Kimura model. The binary model is very well understood and has a lot of special properties \cite{BW}. In particular from \cite[Theorem 4.9]{DBM} we know that Conjecture \ref{glhip} holds for $G=\z_2$.

We have got three natural projections $f_i:\z_2\times\z_2\rightarrow \z_2$ for $i=1,\dots,3$. The map $f_1\times f_2\times f_3:\z_2\times\z_2\rightarrow \z_2\times\z_2\times\z_2$ is injective. Moreover it induces a dominant map from the product of three binary models onto the 3-Kimura model. This map is the key tool that will allow us to transfer some of the properties from the binary model to the 3-Kimura model. Unfortunately the map is not surjective, but just dominant. We can projectivise the varieties, but then we get a rational map. It turns out that a combine use of both of the maps allows to derive the main theorem.

Let $f_i^*:M_S^{\z_2\times \z_2}\rightarrow M_S^{\z_2}$ be a morphism of lattices induced by $f_i$. Let $g:M_S^{\z_2}\rightarrow M_E^{\z_2}$ be the morphism of lattices that corresponds to the parametrization map of the binary model.

We have got the following commutative diagram similar to Diagram \ref{diag}:
$$
\xymatrix{
M_{S}^{\z_2}\times M_{S}^{\z_2}\times M_{S}^{\z_2}\ar[rr]^{g\times g\times g}&&M_E^{\z_2}\times M_E^{\z_2}\times M_E^{\z_2}\\
M_{S}^{\z_2\times \z_2}\ar[rr]^{g_0}\ar[u]^{f_1^*\times f_2^*\times f_3^*}&&M_E^{\z_2\times \z_2}\ar[u]^i\\
}
$$
The following Fact is well known.
\fa\label{dimension3Kbin}\emph{
The dimension of the affine 3-Kimura model is equal to $3|E|+1$.
The dimension of the product of three affine binary models is equal to $3(|E|+1)$.
The dimension of the projective 3-Kimura model is equal to $3|E|$.
The dimension of the product of three projective binary models is equal to $3|E|$.
$\hfill\square$}
\kfa
It follows that if we consider projective varieties representing the models, the dominant morphism described above becomes a rational finite map.
As the map between projective varieties is not a morphism we will restrict our attention only to dense orbits of the tori. On these tori orbits all maps are well defined and are represented by morphism of lattices.
\subsection{Maps of dense torus orbits}
Let us consider the following diagram:\begin{equation}\label{diagproj}\end{equation}

\resizebox{\textwidth}{!}{
\xymatrixcolsep{0pt}
\xymatrix{
\small{M_{S}^{\z_2}\times M_{S}^{\z_2}\times M_{S}^{\z_2}}\ar[rrr]^{g\times g\times g}&&&M_E^{\z_2}\times M_E^{\z_2}\times M_E^{\z_2}&\\
&M_{S,0}^{\z_2}\times M_{S,0}^{\z_2}\times M_{S,0}^{\z_2}\ar@{_{(}->}[lu]\ar[rrr]&&&\widehat M_{E,0}^{\z_2}\times\widehat M_{E,0}^{\z_2}\times \widehat M_{E,0}^{\z_2}\ar@{_{(}->}[lu]\\
M_{S}^{\z_2\times \z_2}\ar[rrr]^(0.55){g_0}\ar[uu]^{f_1^*\times f_2^*\times f_3^*}&&&M_E^{\z_2\times \z_2}\ar@{^{(}->}[uu]^(.43)i&\\
&M_{S,0}^{\z_2\times\z_2}\ar@{_{(}->}[lu]\ar[uu]^(.43)f\ar[rrr]^h&&&\widehat M_{E,0}^{\z_2\times\z_2}\ar@{_{(}->}[lu]\ar@{^{(}->}[uu]^(.55)j\\
}
}

First let us explain the morphism $j$. It is injective, as it is a restriction of $i$. 
The lattice $\widehat M_{E,0}$ is the character lattice of the torus acting on the projective toric variety representing the model. The morphism $j$ is induced by the rational finite map form the product of three $\p(X(T,\z_2))$ to $\p(X(T,\z_2\times\z_2))$. Due to the coordinate system we can identify dense torus orbits with the tori.
\dfi[The torus $\mathbb{T}_X$]
Let $X$ be any toric variety in an affine or projective space with a distinguished coordinate system. We may identify the dense torus orbit of $X$ with the torus using the coordinate system. We denote this torus $\mathbb{T}_X\subset X$.
\kdfi
\fa\label{torus}\emph{
The torus $\mathbb{T}_X$ contains a point of an affine or projective space if and only if this is a point of $X$ with all coordinates different from $0$.}
\kfa

The morphism $j$ of character lattices is induced by the finite morphism from $\mathbb{T}_{(\p(X(T,\z_2)))^3}=(\mathbb{T}_{\p(X(T,\z_2))})^3$ to $\mathbb{T}_{\p(X(T,\z_2\times\z_2))}$. Due to the discussion in the proof of Proposition \ref{inkluzja} we also know that the morphism of ambient spaces does not depend on the tree, but only on the number of leaves $l$. Hence the vertical morphisms of lattices on the left hand side of Diagram \ref{diagproj} are the same for all trees with $l$ leaves.

\subsection{Idea of the proof}\label{idea}
The main reason for passing to tori is that we want to have a well defined dominant finite map. This allows us to take advantage of toric geometry. For example we know that the number of points in the fiber of the morphism of tori is equal to the index $I_1$ of the image of $j$ in $(\widehat M_{E,0}^{\z_2})^3$.

For the projective ambient spaces the situation is a little bit different. The morphism $f$ is not injective, so the corresponding morphism of tori is not surjective. We will show that the image of $f$ in $(M_{S,0}^{\z_2})^3$ is of finite index, say $I_2$. It means that the corresponding morphism of tori is finite with each fiber having $I_2$ elements. Moreover we will show that $I_2=I_1$. Hence we get the diagram:
\label{diagtor}
$$
\xymatrix{
\mathbb{T}_{(\p(W_E^{\z_2})^3)}\ar[r]&\mathbb{T}_{\p(W_E^{\z_2\times\z_2})}\\
\mathbb{T}_{(\p(X(T,\z_2)))^3}\ar@{^{(}->}[u]\ar@{->>}[r]&\mathbb{T}_{\p(X(T,\z_2\times\z_2))}\ar@{^{(}->}[u],\\
}
$$
where the horizontal maps are finite, \'etale of the same degree.

This means that if we consider the morphism of projective ambient spaces, then the preimage of $\mathbb{T}_{\p(X(T,\z_2\times\z_2))}$ is precisely $\mathbb{T}_{(\p(X(T,\z_2)))^3}$. Hence any intersection results that hold for the binary model must also hold for the $3$-Kimura model.
In particular as Conjecture \ref{glhip} holds for the binary model we obtain a set-theoretic version of Conjecture \ref{drhip} for the $3$-Kimura model.

By easy algebraic arguments we will also prove Conjecture \ref{drhip} scheme-theoretically for $3$-Kimura model.
\subsection{Proof}
Our first step will be to understand the morphism of projective ambient spaces $(\p(\widetilde W_E^{\z_2}))^3\dashrightarrow \p(\widetilde W_E^{\z_2\times \z_2})$. This is a well defined map on dense tori orbits. The map of tori corresponds to morphism of lattices $f:M_{S,0}^{\z_2\times\z_2}\rightarrow (M_{S,0}^{\z_2})^3$. This morphisms depend only on the number of leaves, not on the tree.

Each socket is an association of a group element to a leaf. Hence we can embed the group of sockets $\mathfrak{S}$ in $G^l$. We can also view the group $\mathfrak{S}$ as a $\z$-module. This gives us group morphisms $M_S\rightarrow\mathfrak{S}\rightarrow G^l$. The element of the basis of $M_S$ indexed by a socket $s$ is mapped to the socket $s$.
\ex[The case of the binary model and trivalent claw tree]
\emph{
Let us consider the tree $K_{1,3}$ and the group $\z_2$. We have got $4$ sockets: $(0,0,0),(1,1,0),(1,0,1),(0,1,1)$. By coordinate-wise action they form a subgroup of $(\z_2)^3$. The lattice $M_S$ is freely generated by four basis vectors $e_{(0,0,0)},$ $e_{(1,1,0)},$ $e_{(1,0,1)},$ $e_{(0,1,1)}$. The morphism $M_S\rightarrow\mathfrak{S}$ maps $e_{(a,b,c)}$ to $(a,b,c)$. Of course $ke_{(a,b,c)}$ is mapped to $k(a,b,c)$. For example $3e_{(1,1,0)}$ is mapped to $(1,1,0)+(1,1,0)+(1,1,0)=(1,1,0)$.}
\kex
\lem\label{dokladny}
We have an exact sequence of groups:
$$M_{S,0}^{\z_2\times\z_2}\rightarrow (M_{S,0}^{\z_2})^3\rightarrow (\z_2)^l.$$
The first morphism is given by $f$. The second is the group operation on sockets described above\footnote{In this case the second operation is often called XOR.}.
\klem
\dow
It is clear that this is a complex. Let $(b_i')_{i\geq 0}$ be the basis of $M_{S}^{\z_2}$ corresponding to sockets. Let $s_i$ be the socket corresponding to $b_i'$. Moreover suppose that $b_0'$ corresponds to the trivial socket, that is the neutral element of $\mathfrak{S}$. Let $b_i$ be the basis of $M_{S,0}^{\z_2}$ defined as $b_i=b_i'-b_0'$ for $i>0$. Note that an element $(b_i',b_j',b_k')$ is in the image of $f_1^*\times f_2^*\times f_3^*$ iff the corresponding three sockets $s_i$, $s_j$, $s_k$ sum up to the neutral element of $\mathfrak{S}$. Hence the elements of the form $(b_i,b_i,0)=(b_i',b_i',b_0')-(b_0',b_0',b_0')$ are in the image of $f$. We see that $(2b_i,0,0)=(b_i,b_i,0)+(b_i,0,b_i)-(0,b_i,b_i)$ is also in the image. Furthermore for any two sockets $s_i$ and $s_j$ there exists a socket $s_k:=s_i+s_j$ such that $(b_i,b_j,b_k)$ is in the image of $f$. This reduces any element from $(M_{S,0}^{\z_2})^3$ to an element $(b_i,0,0)$ modulo the image of $f$ or to $0$. Hence any element is in the image if the XOR of all its coordinates is zero.
\kdow
\dfi[The kernel $K$]
Let us consider the restriction $K=K_1\times K_2\times K_3\subset M_{S,0}^{\z_2}\times M_{S,0}^{\z_2}\times M_{S,0}^{\z_2}$ of the kernel of the morphism
$g\times g\times g$.
\kdfi
Each character in $K$ is a character of $(\mathbb{T}_{\p(\widetilde W_E^{\z_2})})^3$, that is the trivial character when
restricted to the product $(\mathbb{T}_{\p(X(\z_2))})^3$. Each such character is a triple of characters of $\mathbb{T}_{\p(\widetilde W_E^{\z_2})}$. Each character of the triple is a quotient of monomials $\frac{m_1}{m_2}$ of the same degree on the projective space $\p(\widetilde W_E^{\z_2})$.
The polynomials $m_1-m_2$ span\footnote{They do not only generate the ideal, but even span it as the vector space.}
 the ideal of the toric variety $\p(X(\z_2))$. We want to view characters as functions. Hence we restrict our attention to $(\mathbb{T}_{\p(\widetilde W_E^{\z_2})})^3$. In the algebra of this torus the ideal of $(\mathbb{T}_{\p(X(\z_2))})^3$ is generated by elements $k-1$, where $k\in K$.
\dfi[The kernel $D$]
For any tree $T$ let $D^T$ be the kernel of the map $h$.
\kdfi
The elements of $D$ represent characters trivial on the projective 3-Kimura variety. In the setting of Subsection \ref{Intoftor} we want to prove that sublattices $D^T$ for different trees $T$ with $l$ leaves generate the sublattice $D^{K_{1,l}}$. The idea is to push the lattices $D$ to $(M_{S,0}^{\z_2})^3$ using the morphism $f$. Next we use the results on binary models to obtain the generation for $f(D)$. Using properties of the image of $f$ we are able to conclude the generation in $M_{S,0}^{\z_2\times\z_2ć}$. The following lemma enables us to restrict to the image of $f$ instead of regarding whole lattice $(M_{S,0}^{\z_2})^3$. Algebraically it means that characters of $\mathbb{T}_{\p(\widetilde W_E^{\z_2\times\z_2})}$ trivial on the image of $(\mathbb{T}_{\p(\widetilde W_E^{\z^2})})^3$ are also trivial on $\mathbb{T}_{\p(X(T,\z_2\times\z_2))}$ for any tree $T$ with $l$ leaves. Hence geometrically the image of $(\mathbb{T}_{\p(\widetilde W_E^{\z^2})})^3$ contains $\mathbb{T}_{\p(X(T,\z_2\times\z_2))}$.
\lem\label{gl}
For any tree $T$ the kernel $K^T$ is a sublattice of the image of $f$.
\klem
\dow
It is enough to show that $K_1\times\{0\}\times\{0\}\subset \Imi f$. Suppose that $m=\sum_i a_i b_i\in K_1$, where each $b_i$ is as in the proof of Lemma \ref{dokladny}. Hence $b_i=(g_1^i-e,\dots,g_l^i-e)$, where $e$ is the neutral element of $\z_2$ and $g^i_j\in \z_2$ are elements forming a socket. We know that $g(m)=0$. In particular the coordinates of $M_E$ indexed by leaves are equal to zero. Let us fix $k$ that is a number of a leaf $1\leq k\leq l$.
Let us look at all coordinates indexed by pairs $(k,q)$ where $q\in \z_2$. The restriction of $M_E$ to these coordinates is a free abelian group spanned by elements of $\z_2$. Hence $\sum_i a_i(g_k-e)=0$ in the free abelian group generated formally by elements of $\z_2$. Hence, a fortiori, $\sum_i a_i(g_k-e)=e$ where now the sum is taken in $\z_2$. As the action in $\mathfrak{S}$ is coordinate-wise we see that the image of $m$ in $\mathfrak{S}$, and hence in $\z_2^l$, is the neutral element. Using Lemma \ref{dokladny} we see that $m\in \Imi f$.
\kdow

\ob
The index of the image of $f$ in $(M_{S,0}^{\z_2})^3$ is equal to  the index of the image of $g$ in $(\widehat M_{E,0}^{\z_2})^3$.
\kob
\dow
This is a consequence of the injectivity of $g$ and Lemma \ref{gl}.
\kdow
\wn
Conjecture \ref{drhip} holds set-theoretically.
\kwn
\dow
The index of the image of of $f$ equals the degree of the finite map of tori. In particular we are in the situation of Diagram \ref{diagtor}. The corollary follows from the discussion at the beginning of Section \ref{idea}.
\kdow
Now we will prove Conjecture \ref{drhip} scheme-theoretically. Let $T_0=K_{1,l}$. We consider such trees $T_i$ that the ideal of $\mathbb{T}_{\p(X(T,\z_2))}$ is the sum of the ideals $\mathbb{T}_{\p(X(T_i,\z_2))}$. Let $K^{T_i}$ be the kernel of $g\times g\times g$ for the tree $T_i$. Let $D^{T_i}$ be the kernel of $h$ for the tree $T_i$. We know form Subsection \ref{Intoftor} that the lattices $K^{T_i}$ for $i>0$ span $K^{T_0}$.
\tw
The lattices $K^{\z_2\times\z_2}_i$ for $i>0$ span $K^{\z_2\times\z_2}_0$. Conjecture \ref{drhip} holds scheme theoretically.
\ktw
\dow
Let $a\in K^{\z_2\times\z_2}_0$. We know that $f(a)\in K^{\z_2}_0$, so $f(a)=\sum k_i$, where $k_i\in K^{\z_2}_i$. Using Lemma \ref{gl} we can find $k_i'\in K^{\z_2\times\z_2}_i$ such that $f(k_i')=k_i$. This means that $a-\sum k_i'$ is in the kernel of $f$. In particular, as $g$ is injective, $a-\sum k_i'$ belongs to every $K^{\z_2\times\z_2}_i$, hence we obtain the desired decomposition.
\kdow
\uwa
From Appendix \ref{A} it is enough to take two (particular) different $i>0$ to span $K^{\z_2\times\z_2}_0$, as it was in the case of binary model.
\kuwa

\section{Applications}\label{Ap}
In this section we present a few applications of the main theorem. The basic result that we use is due to Marta Casanellas and Jes\'us Fern\'andez-S\'anchez \cite{CFS}. It states that all biologically meaningful points are contained in the dense torus orbit of $X(T,\z_2\times\z_2)$. In the main theorem we give a precise description of this orbit for any tree. This is sufficient for biologists.

People dealing with applications are usually interested in trivalent trees. Let us motivate the use of other trees. The first, obvious reason is that they can appear (at least hypothetically) as right models of evolution. This however is a degenerate situation that is often neglected. The next subsection presents a different reason.
\subsection{Identifiability}
 Dealing with applications we are given a point $P$ in the space of all possible probabilities $\widetilde W_L$. The first question is for which trees this point can be realized. More precisely for which trees $T$ we have an inclusion $P\in X(T,\z_2\times\z_2)$. We are interested in knowing if this is only one tree $T$ or there are several possibilities. This is a first part of the \emph{identifiability} problem. Hence Conjecture \ref{glhip} is a question about the locus of points for which the identifiability problem cannot be resolved at all. Of course a generic point that belongs to any of the varieties belongs to exactly one $X(T,\z_2\times\z_2)$ with $T$ trivalent. Much more is known about the identifiability of different models. For the precise results the reader is advised to look in \cite{Theidentoftreetopol} or \cite{Identof2tree} and the references therein.

In particular we see that points that belong to some $X(T,\z_2\times\z_2)$ where $T$ is not trivalent cannot identify the tree topology. Hence the question about the locus of these points, or equivalently about the polynomials defining such varieties may give some results for trivalent trees. However, as situation in Appendix \ref{B} shows, the phylogenetic invariants of two varieties $X(T,\z_2)$ for two different trees, do not generate the ideal of the variety associated to their degeneration.

The second, but equally important question about the identifiability is to give the description of the fiber $\psi^{-1}(P)$. The biologist aim at distinguishing one point in the fiber. This would enable to identify not only the tree topology, but also corresponding probabilities of mutation. The algebraic setting allows us to give a description of this fiber. We assume that $P$ is biologically meaningful, that is is contained in the dense torus orbit. Equivalently all coordinates of $P$ after the Fourier transform are different from zero. We prefer to work up to multiplicity, that is regard the projectivization of $\psi$ denoted by $\psi_\p$. The fiber $\psi_\p^{-1}(P)$ is contained in the dense torus orbit of $\prod \p(W_e)$. As this parameter space is of the same dimension as image, we know that $\psi_\p$ is a generically finite map. Moreover when restricted to dense torus orbits it is \'etale and finite. Hence each fiber is finite and contains the same number of points, independent from $P$. This number is the index of lattice $\widehat M_E$ in a saturated sublattice of $M_E$. For 3-Kimura model this number can be calculated from \cite{CFS}. Of course we do not claim that all the points in the fiber have got a probabilistic meaning. We just prove that from the algebraic point of view there is always a fixed, finite number of possible candidates for transition matrices.
\subsection{Phylogenetic invariants}
The main theorem gives an inductive way of obtaining phylogenetic invariants of any tree. It is an open problem if these invariants generate the whole ideal. It is proved however that they give a description of all biologically meaningful points.
The method is very simple. Suppose that we know the phylogenetic invariants for all trees with vertices of degree less or equal to $d$. Due to the results of \cite{SS} it is enough to describe the phylogenetic invariants for the claw tree $K_{1,d+1}$. To obtain them we just take the sum of two ideals. They are both associated to trees with the same topology. The tree has got two inner vertices $v_1$ and $v_2$ of degrees 3 and $d$ respectively. The difference between the ideals is a consequence of different labeling of leaves. For one tree the leaves adjacent to $v_1$ are labelled by 1 and 2. For the second tree they are labeled 1 and 3. Notice that in fact we have to compute just one ideal. The second one can be obtained by permuting the variables.
\section{Special cases}\label{Sc}
\subsection{Orbits}
This section contains technical results concerning some orbits of the torus action that are not dense in the variety. Such orbits correspond to faces of the polytope associated to the variety. Each orbit is isomorphic to a torus. Its character lattice is the lattice spanned by points of the corresponding face of the polytope. In this notation the dense orbit corresponds to the polytope.

Using the coordinates this description can be made more explicit. The points of the polytope correspond to coordinates of the affine (or projective) space. By choosing a face $F$ we distinguish coordinates corresponding to points of this face. The orbit corresponding to $F$ is the projection of the dense orbit onto the linear subspace. This subspace is given by a condition that all coordinates corresponding to points that are not in $F$ are equal to zero. Hence closures of all orbits are given as intersections of the variety with linear subspaces corresponding to faces. The orbits themselves correspond to projections of the dense orbit onto these subspaces. In particular let $Q$ be a point with all coordinates equal either to $0$ or $1$. This point is contained in the toric variety if and only if points of the polytope corresponding to coordinates with the entry equal to $1$ form a face. In case $Q$ is contained in the toric variety it is the unique distinguished point of the orbit. It corresponds to the neutral element of the torus after identifying it with the orbit.

We can decompose the affine space $\a^n$ into tori $S_i=\{0\}^k\times {\c^*}^{n-k}$ by setting some coordinates equal to zero and requiring that others are different form zero. A toric variety intersected with such a subset is either empty or contains the unique distinguished point $Q_i=\{0\}^k\times\{1\}^{n-k}$. If it is nonempty it is a torus. If we intersect a few toric varieties with $S_i$ the situation is a little bit more complicated. We can get an empty set if and only if one of the varieties has an empty intersection with $S_i$. Otherwise the intersection contains the point $Q_i$. It might however be reducible.

To be even more precise let $M$ be the character lattice of the torus acting on $\a^n$. Let $M_i$ be the character lattice of affine toric varieties $X_i\subset \a^n$. Let $\varphi_i:M\rightarrow M_i$ be the surjective morphism corresponding to the inclusion and let $K_i$ be its kernel. There is a natural projection of $\a^n\supset (\c^*)^n\rightarrow S_i$. This induces an inclusion of lattices $M'\subset M$, where $M'$ is the character lattice of the torus acting on $S_i$. We may restrict $\varphi_i$ and $K_i$ to $M'$ obtaining $\varphi_i'$ and $K_i'$. The question whether $\bigcap X_i\cap S_i$ is irreducible is equivalent to the question if $\bigcup K_i'$ is saturated in $M$. If not, its index is equal to the number of connected components\footnote{Each is irreducible and isomorphic to a torus.}. Moreover the saturation of $K_0:=\bigcup K_i'$ gives characters trivial on the distinguished component of the intersection that contains $Q_i$. Hence $M/K_0$ is the character lattice of this component.

To check Conjecture \ref{glhip} set theoretically we may check it separately on each $S_i$. This approach is quite standard, see e.g.~\cite{EiS}.
For a given $S_i$ we have a few different possibilities.

Case 1) The intersection is empty.

The first possibility is that the intersection $\bigcap_{T\neq K_{1,l}}X(T,G)\cap S_i$ is empty. It follows that $X(T,G)\cap S_i$ is empty for some $T\neq K_{1,l}$. As $X(K_{1,n},G)\subset X(T,G)$ the intersection $X(K_{1,n},G)\cap S_i$ is also empty hence Conjecture \ref{glhip} holds on $S_i$.

Case 2) The intersection is reducible.

In order to prove the Conjecture \ref{glhip} one has to exclude this possibility.

Case 3) The intersection is nonempty and irreducible.

Here we have to exclude the case of the strict inclusion $X(K_{1,n},G)\cap S_i\subset\bigcap_{T\neq K_{1,l}}X(T,G)\cap S_i$. As both intersections are isomorphic to tori it is enough to compare their dimensions.
\subsubsection{Special faces - type 1}
Unfortunately we are able do conclude only for some special tori $S_i$. Moreover the reasoning presented in this section is often very technical. We would much appreciate a uniform approach that would enable the proof of Conjecture \ref{glhip}.

First let us consider faces $F_{(e,g)}$ of the type $h_{(e,g)}^*(x)=1$, where $h_{(e,g)}\in M_E$ is a basis vector. This face contains networks that associate $g$ to edge $e$. Let us consider the claw tree $K_{1,l}$ and any relation $R$ between the vertices of the face $F_{(e,g)}^{K_{1,l}}$. Let $T$ be a tree with two inner vertices. One of the vertices $v_1$ is of degree $3$ and the second one is of degree $l-1$. The edge corresponding to $e$ is a leaf adjacent to $v_1$:
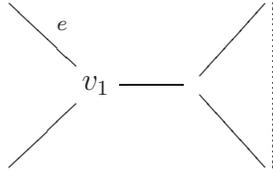
\begin{equation}\label{drzewo}\end{equation}
$$\xymatrix{
\ar@{-}^e[dr]&&&\ar@{-}[dl]\ar@{..}[dd]\\
&v_1\ar@{-}[r]&\ar@{-}[dr]&\\
\ar@{-}[ur]&&&\\
}$$
One can easily see that the relation $R$ lifts up to the relation between the vertices of the polytope associated to $T$. Hence $X(K_{1,l},G)$ and $X(T,G)$ are scheme theoretically equal on the closure of $S_i$. Faces contained in some $F_{(e,g)}$ will be called faces of type 1.
\subsubsection{Special faces - type 2}
Now let us assume that $G=\z_2\times\z_2$. We will give a description of faces that contain a diagonal. We say that a face contains a diagonal if for some two vertices, it is a minimal face containing them. Hence all edges are included. Examples of polytopes that do not contain a diagonal are simplices or piramids.

Let us consider two networks $n_1$ and $n_2$ on a claw tree. We will use the same notation for networks and corresponding vertices of the polytope.
\ob
The minimal face containing $n_1$ and $n_2$ contains precisely those networks $n$, such that for each edge $e$ we have either $n(e)=n_1(e)$ or $n(e)=n_2(e)$.
\kob
\dow
It is obvious that this set of networks forms a face. If we consider any network $n$ with the property above, we may define a new network $n'$ by the property:
$$
n'(e)=\begin{cases}
n_1(e)\text{  if  }n(e)=n_2(e)\\
n_2(e)\text{  if  }n(e)=n_1(e)\\
\end{cases}.
$$
For corresponding vertices we have $n+{n'}={n_1}+{n_2}$. Hence $n$ and ${n'}$ must belong to to the minimal face containing $n_1$ and $n_2$.
\kdow
For $S_i$ corresponding to such faces one can also prove Conjecture \ref{glhip}. The proof is technical and amounts to decomposing any network in the face plus $kn_1$ for $k\in\n$ into sum of networks that do not differ much from $n_1$.
\subsection{Irreducibility}
We finish with a remark about the reducibility of the intersection of varieties associated to different trees. We prove that for the 3-Kimura model if the intersection is reducible, then it must have at least three components.

The space $\widetilde W_E$ is the regular representation of $\mathfrak{N}$. The variety $X(T,G)$ is invariant with respect to its action. We also know
that each component of the intersection $\bigcap X(T,\z_2\times\z_2)$ different from $X(K_{1,l},\z_2\times\z_2)$ is in a hyperplane section.
Suppose that Conjecture \ref{glhip} does not hold set theoretically. 
The action of $\mathfrak{N}$ permutes the components. It also acts transitively on variables. In particular it sends a hyperplane section to a different hyperplane section. Hence there must be at least several components in the intersection.
\section{Appendix}
\subsection{}\label{A}
Let us fix a number of leaves $l$. We claim that for special two trees $T_1$ and $T_2$ the scheme-theoretic intersection $X(T_1,\z_2)\cap X(T_2,\z_2)$ equals $X(K_{1,l},\z_2)$. We number the leaves from $1$ to $l$. The trees $T_1$ and $T_2$ are isomorphic as graphs but have different leaf labeling. The topology of the trees is the same as the one described in Picture \ref{drzewo}. For the tree $T_1$ the leaves adjacent to $v_1$ have got numbers $1$ and $2$. For the tree $T_2$ they are numbered $1$ and $3$. The ideal of the variety associated to a tree for the group $\z_2$ is always generated in degree 2 \cite{Sonja}, \cite{SS}. Hence the generators of the ideals are of the form $n_1n_2=n_3n_4$ where $n_i$ for $1\leq i\leq 4$ are coordinates corresponding to networks. A network is an association of group elements to edges. In can be represented by a vector, whose entries are group elements, indexed by edges. A product of coordinates corresponding to networks can be represented by a matrix with columns representing the networks. Of course we can permute the columns in such a matrix. An equality of two products can be represented by two matrices. In fact each binomial equality corresponds to a pair of matrices $(M_0,M_1)$ whose columns represent networks and rows are the same up to permutation. Hence each generator of the ideal of $X(K_{1,l},\z_2)$ is represented by a pair of $2\times l$ matrices with entries from $\z_2$. Moreover the sum in each column is the neutral element and rows of both matrices are the same up to permutation. As we can permute columns of each matrix we may assume that the first rows of both matrices coincide. Let us consider any such generator $(M_0,M_1)$ in the ideal of $X(K_{1,l},\z_2)$.

First suppose that the entries in the first row are the same, that is either $00$ or $11$. Then the relation holds both for $X(T_1)$ and $X(T_2)$. Hence we may suppose that the first row is $01$ or $10$. The same reasoning holds for the second and third row. Hence all three rows in both matrices are either $01$ or $10$. If the second (resp.~third) rows are the same in both matrices then the relation holds for $X(T_1)$ (resp.~$X(T_2)$). So the only possibility left is that the second and third rows of $M_1$ are the negation of the second and third rows of $M_0$. In this case the relation does not hold in any $X(T_i)$ but we can generate it. We consider a matrix $M$ that is equal to $M_0$ with the first two rows permuted. The pair $(M_0,M)$ represents a relation in $X(T_1)$. Moreover the pair $(M,M_1)$ represents a relation in $X(T_2)$.
\subsection{}\label{B}
We will prove that a scheme theoretic intersection $X(T_1,\z_2)\cap X(T_3,\z_2)$ does not have to be equal to $X(K_{1,l},\z_2)$ even if $K_{1,l}$ is the only tree smaller then $T_1$ and $T_3$. We consider the case of five leaves $l=5$. The tree $T_1$ is the same as in the previous subsection. The tree $T_3$ is isomorphic, with two distinguished leaves labeled with $4$ and $5$. We consider the relation given by a pair of matrices:
\[
\left[
\begin{array}{cccccccc}
1&0\\
0&1\\
0&0\\
0&1\\
1&0\\
\end{array}
\right],
\left[
\begin{array}{cccccccc}
1&0\\
1&0\\
0&0\\
0&1\\
0&1\\
\end{array}
\right].
\]
This corresponds to a generator of the ideal of $X(K_{1,5},\z_2)$. Consider any relation involving the first matrix and some other matrix $M$ for $X(T_1)$ or $X(T_3)$. One can see that the first two rows of $M$ must be negations of each other and the third one is $00$. Hence it is impossible to generate the relation above.

\bibliographystyle{amsalpha}
\bibliography{xbib}
{
\footnotesize
\noindent
Mateusz Micha\l ek\nopagebreak

\noindent
{Mathematical Institute of the Polish Academy of Sciences,\nopagebreak

\noindent
\'{S}niadeckich 8, 00-956 Warszawa, Poland}
\vskip 2pt
\noindent
{Institut Fourier, Universite Joseph Fourier,

\noindent
100 rue des Maths, BP 74, 38402 St Martin d'H\`eres, France}

\noindent
e-mail address:\emph{wajcha2@poczta.onet.pl}}
\end{document}